\documentclass[]{interact}

\usepackage{epstopdf}
\usepackage[caption=false]{subfig}

\usepackage{mathrsfs}
\usepackage{amsmath}
\usepackage[font=small,labelfont=bf,labelsep=none]{caption}
\usepackage{threeparttable}
\usepackage{tabularx}
\usepackage{subfig}
\usepackage{booktabs}
\usepackage{array}
\usepackage{rotating}

\usepackage{varwidth}
\usepackage{graphicx}

\usepackage{makecell}
\usepackage{multirow}
\usepackage{pdflscape}

\usepackage[numbers,sort&compress]{natbib}
\bibpunct[, ]{[}{]}{,}{n}{,}{,}

\theoremstyle{plain}

\theoremstyle{definition}

\theoremstyle{remark}
\newtheorem{remark}{Remark}

\renewcommand{\figurename}{Fig.}

\begin{document}

\articletype{ARTICLE TEMPLATE}


\title{Globally Coupled Particle Swarm Optimization}

\author{
\name{Liguo Yuan\textsuperscript{a}\thanks{CONTACT L.~G. Yuan. Email: liguoy@scau.edu.cn}}
\affil{\textsuperscript{a}Department of Mathematics, College of Mathematics and Informatics, South China Agricultural University, Guangzhou, China.}
}

\maketitle

\begin{abstract}
All things in the world are interconnected, the only difference is the strength of their connections.
Particle swarm optimization(PSO) simulates the foraging behavior of a flock of birds, information is transmitted to quickly find the location of food.
This is a process of information exchange, where birds influence each other, constantly adjusting their position and speed values,
and updating the optimal position information.
In this paper, we propose a globally coupled particle swarm optimization(GCPSO) that combines the globally
coupled map lattices(GCML) with the PSO to enhance its optimization capabilities.
The information of the $i$-th lattice point is influenced by the information of all lattice points in GCML.The information between lattice points is interdependent.
Inspired by this, we will integrate GCML into PSO and propose a new improved particle swarm optimization, namely GCPSO.
Here, the speed update formula has been modified and improved.
The next flight speed of each bird is influenced by its current speed, its historical best position, the historical global best position, its current position, and the positions of all  birds.
Each bird(particle) is influenced by all birds. The strength of the impact will be distinguished by the size of the weight. This is the essential difference and key improvement from PSO.
Through extensive experiments, it has been found that compared to PSO, GCPSO has a stronger ability to search for solutions.
\end{abstract}


\section{Introduction}\label{sec_intro}
Optimization problem is the process of finding the optimal solution, with regard to a given criterion, from a set of available solutions.
The global optimization problem is usually defined as
\begin{equation}
X^*=\arg \min\limits_{X\in S} f(X), \label{miniproblem}\\
\end{equation}
where
$X=(x_1, x_2, \cdots, x_n)$ and  $S=[a_1, b_1] \otimes  [a_2, b_2] \otimes  \cdots  \otimes [a_n, b_n].$
The objective function $f(X)$ is assumed to be continuous and differentiable \cite{parallelPSO,30objectivefunctions}. Optimization algorithms require gradient  information which is depicted as
$\nabla f(X)=[\frac{\partial f}{\partial x_1}(X),\frac{\partial f}{\partial x_2}(X), \cdots, \frac{\partial f}{\partial x_n}(X)]$, then the differentiability of $f(X)$ is a necessary condition.
However, swarm intelligence optimization algorithms(such as PSO \cite{Kennedy1995,Eberhast1995})  break these limitations.
$f(X)$ can be non-differentiable or even discontinuous.
Optimization problems arise in various quantitative disciplines. There is a wide range of problems in the relevant literature that can be treated as global optimization problems
and these problems can be successfully solved by PSO.
Anula Khare et al. reviewed the PSO and successfully applied it in the optimization problem of Solar Photovoltaic system \cite{Anula}.
Houssein et al. conducted a rigorous and systematic comprehensive analysis of the PSO and including the diverse applications of the algorithm \cite{Houssein}.
AlRashidi et al. presented a comprehensive coverage of different PSO applications in solving optimization problems in the area of electric power systems.
and discusseed PSO possible future applications in the area of electric power systems and its potential theoretical studies \cite{AlRashidi}.
Russell C.Eberhart and Yuhui Shi analyzed the engineering and computer science aspects of developments, applications, and resources related to PSO \cite{Shi}.
Mr. Ninad K. Kulkarni et al.  reviewed the applications of PSO  in mechanical domain and also described its improved version \cite{Kulkarni}.
Ekrem \"{o}zge and Bekir Aksoy analyzed the trajectory planning of the robotic arm  using the  PSO \cite{Ekrem}.
Ahmed G. Gad systematically reviewed the PSO, including algorithm improvements, diverse application domains, open issues and future perspectives \cite{Ahmed}.
Tareq M. Shami  et al provided a comprehensive review of PSO including remarkable engineering applications \cite{Tareq}.
Tiwari Sukriti and Ashwani Kumar  provided insights on basic notions and progress of the PSO in power system applications \cite{Tiwari}.
Janmenjoy Nayak,et al. presented an in-depth analysis of PSO with its developments from 1995 to 2020 \cite{Nayak}.
John A. Ramirez-Figueroa  et al proposed a new method for disjoint principal component analysis based on PSO \cite{Ramirez}.
Yu Xue et al successfully applied PSO to optimize the topological construction of the  echo state networks \cite{Xue}.
Pace Francesca et al reviewed the application of the PSO to perform stochastic inverse modeling of geophysical data\cite{Pace}.
Ali R. Kashani et al  provided a detailed review of applications of PSO on different geotechnical problems \cite{Kashani}.
There are also many related books that involve theory and application \cite{Felix,Mercang,Micael,Olsson,Mikki,Eberhart}.
PSO has proven to be extremely useful in various disciplines.

On the other hand, there is also a lot of work on improving the PSO. Mainly, PSO has been modified by following strategies:
$(1).$ Dynamic adjustment strategy of inertia weight;
$(2).$ Asynchronous optimization design of learning factors;
$(3).$ Algorithm fusion and multi-strategy collaborative improvement;
$(4).$ Speed and location update mechanism innovation;
$(5).$Parameter adaptation and diversity enhancement, and so on.

Since the introduction of the PSO algorithm in 1995 by Kennedy and Eberhart \cite{Kennedy1995,Eberhast1995}, the
PSO algorithm has attracted a great attention

 modification of the PSO controlling parameters, hybridizing PSO with
other well-known meta-heuristic algorithms such as genetic algorithm (GA) and differential evolution (DE),
cooperation and multi-swarm techniques. This paper attempts to provide a comprehensive review of PSO,
including the basic concepts of PSO, binary PSO, neighborhood topologies in PSO, recent and historical PSO variants,
remarkable engineering applications of PSO, and its drawbacks.

Moreover, this paper reviews recent studies that utilize PSO
to solve feature selection problems. Finally, eight potential research directions that can help researchers further enhance the
performance of PSO are provided.

Davoud Sedighizadeh, et al introduced the Generalized Particle Swarm Optimization algorithm as a new version of the PSO algorithm for continuous space optimization \cite{Sedighizadeh}.


\section{Particle Swarm Optimization}
Particle swarm optimization(PSO) is a swarm intelligent algorithm that simulates the foraging process of animals in nature,
and it is a evolutionary algorithm\cite{PSOKennedy,ShiRussell,Bratton,Sengupta}.
Every candidate solution  is called particle in (\ref{miniproblem}).
All candidate solutions are called swarm.
Each particle is characterized by two factors,
i.e., position $x_i=(x_{i1}, x_{i2}, \cdots, x_{in})$ and velocity $v_i=(v_{i1}, v_{i2}, \cdots, v_{in})$,
where $i$ and $n$ denote the ith particle  and its dimension in search space $\Omega$.
The fitness of each particle can be evaluated according to the objective
function $F(x_i)$ in (\ref{miniproblem}). PSO starts with the random initialization of a swarm in the search space $\Omega$.
Then, every particle update its velocity and position based on its own experience and the experience of optimal particle.
$pbest_i(k)$ represents the optimal position of the $i$-th particle within $k$ iteration steps.
$gbest^k$ denotes the best position in the swarm so far.
The updating rule of velocity and position is shown in (\ref{iteratevelocity}) and (\ref{iterateposition}).
\begin{align}
v_i^{k+1}&=w\times v_i^k+c_1\times r_1\times(pbest_i(k)-x_i^k)+c_2\times r_2\times(gbest^k-x_i^k), \label{iteratevelocity}\\
x_i^{k+1}&=x_i^k+v_i^{k+1},~~~i=1,2,\cdots,n,\label{iterateposition}
\end{align}
where $n$ is the swarm size. $v_i^k$ and $x_i^k$ represent the velocity and position of particle $i$ at $k$-th iteration step respectively.
These mainly are addition and subtraction of vectors, such as $x_i^k=(x_{i1}^k, x_{i2}^k, \cdots, x_{in}^k)$ and $v_i^k=(v_{i1}^k, v_{i2}^k, \cdots, v_{in}^k)$.
$r_1$ and $r_2$ are two independent random numbers between $0$ and $1$.
$c_1$ and $c_2$ are acceleration constants, usually $c_1=c_2=2$, which determine the relative pull of $pbest_i(k)$ and $gbest^k$;
$w$ is called the inertia weight factor which can be chosen by a random number.
Generally, the value of each component in $v_i$ can be clamped to the range $[v_{min}, v_{max}]$.
It controls excessive roaming of particles outside the searching space.
Each particle updates its position according to Eqs.(\ref{iteratevelocity}) and (\ref{iterateposition}).
In this way,  all particles find their new positions and apply these new positions to update their individual best positions and global best position of the swarm.
This process continues until the user specified stopping criterion is met, such as maximum iteration number, the optimal solution tends to stabilize
or total error requirement. The particle $x_i^k$ is drawn towards $pbest_i(k)$ and $gbest^k$.
For more details see references\cite{ PSOKennedy,ShiRussell,Bratton,Poli,Vanneschi,Nayak,Chih,ClercPSO,Shi,overviewpso,comprehensive,Gad,Yuan} and the references cited therein.

In the PSO, each particle calculates its fitness value and compares it with the fitness values of its global neighbors (i.e. the entire particle swarm) to obtain the global optimal value $gbest$ (social influence).
At the same time, each particle calculates its own historical fitness value and compares it with its current fitness value to obtain the optimal fitness value of the individual, that is, the individual extremum $pbest_i$ (learning from
experience).

\begin{figure}[ht]
  \centering
  \subfloat[Update the position]
  {
      \includegraphics[width=0.4\textwidth]{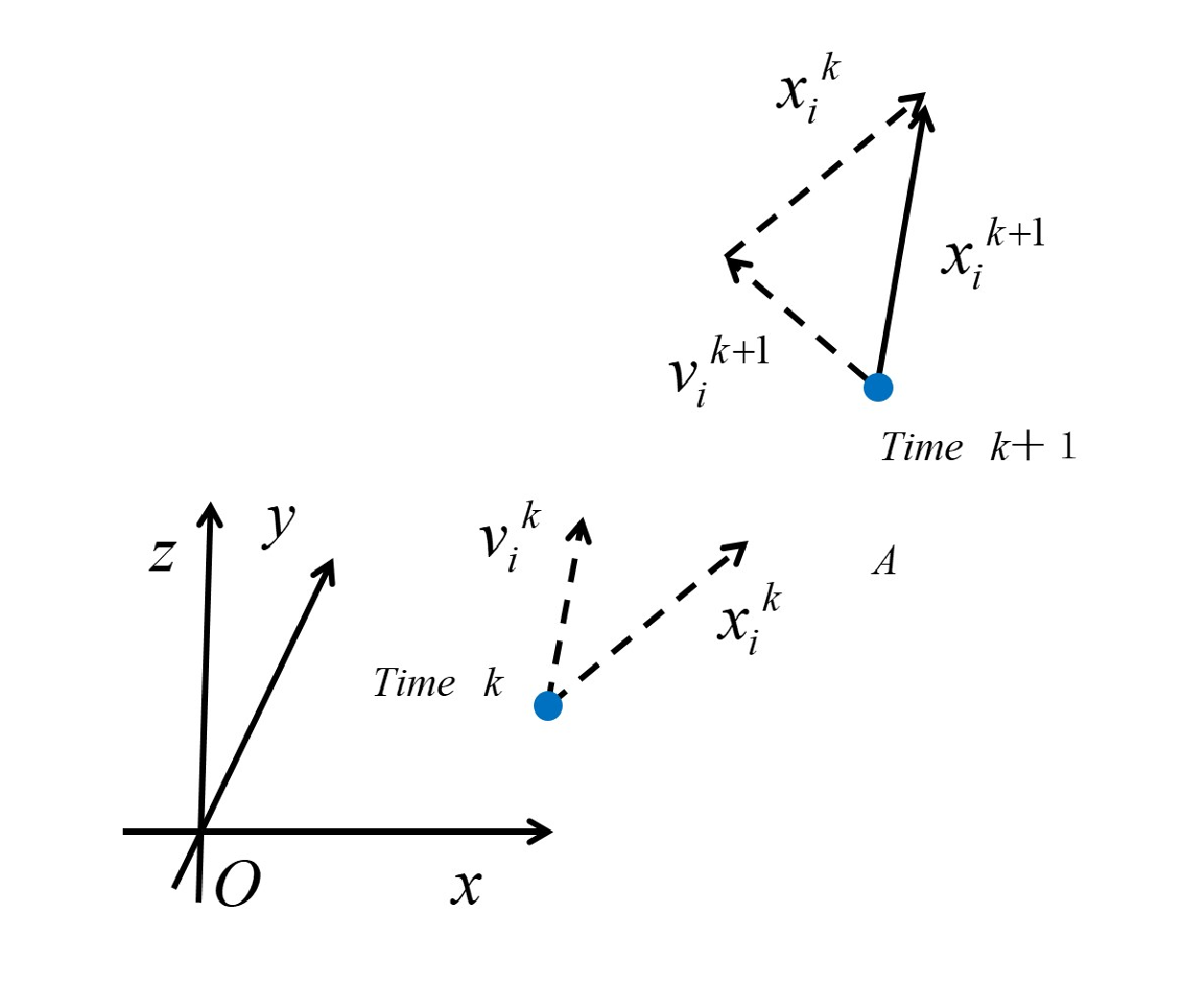}\label{vectoradd}
        }
  \subfloat[Update the velocity]
  {
      \includegraphics[width=0.4\textwidth]{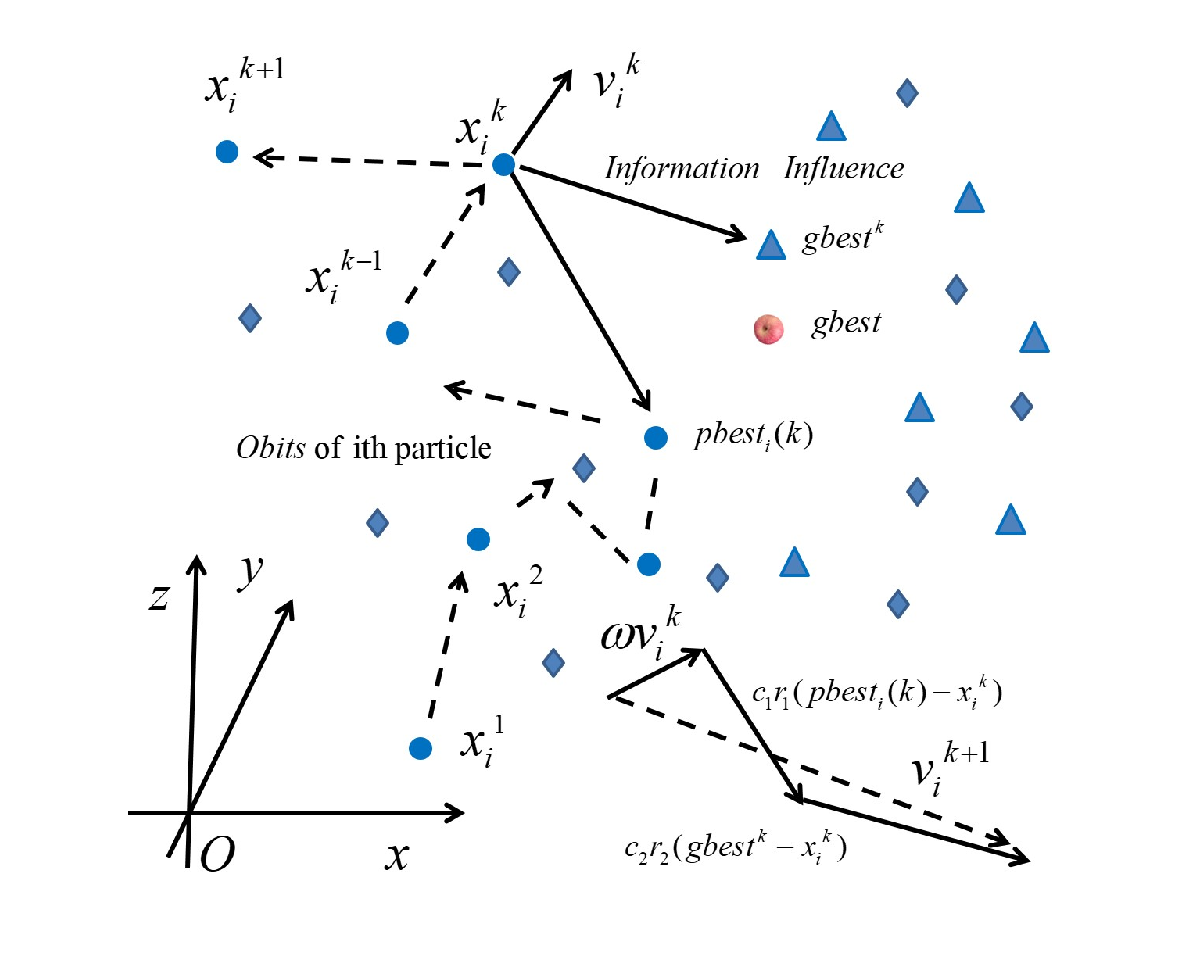}\label{psomovecase}
  }
  \caption{Particle Swarm Optimization  in (\ref{iteratevelocity}) and (\ref{iterateposition})}
  \label{basicpso}
\end{figure}

A vector is a quantity that has two independent properties: magnitude and direction.
Specially, a free vector can move freely in parallel. Both velocity and position are free vectors in $R^n$.
For example, the addition of vectors follows the parallelogram rule (or equivalent to the triangle rule ) in $R^3$, which
is shown in Fig.\ref{basicpso}(a).
The meaning of $pbest_i(k)$, $x_i^k$ and $gbest^k$ of equation (\ref{iteratevelocity}) are shown in Fig.\ref{basicpso}(b).
Food (apple) represents the globally optimal solution. In the previous $k$ iterations, $gbest^k$ represents the current globally optimal solution.
$x_i^1$, $x_i^2$, $\cdots$, $x_i^k$ is the trajectory of the $i$-th particle.
$pbest_i(k)$ is the individual optimal solution of the $i$-th particle.
$pbest_i(k)$ and $gbest^k$ typically possess significant influence over $x_i^k$.
Velocity $v_i^{k+1}$ is influenced by three parts of information in the equation(\ref{iteratevelocity}).
$ w\times v_i^k$ is inertia of velocity. $c_1\times r_1\times(pbest_i(k)-x_i^k)$ represents individual cognition part. $c_2\times r_2\times(gbest^k-x_i^k)$ is social part.
This means that individual $x_i^k$ shares the optimal location information of the swarm.

Equations (\ref{iteratevelocity}) and (\ref{iterateposition}) constitute the main structure of PSO.
The speed of the next generation $v_i^{k+1}$ is influenced by the speed of their previous generation $v_i^k$, its own previous generation position $x_i^k$,
the individual extremum ($pbest_i(k)$) of its previous generation  and the global extremum ($gbest^k$) of the previous generation of the swarm.
These three variables($v_i^k$, $pbest_i(k)$ and $gbest^k$) have the greatest impact on it, and retaining only these three variables simplifies the particle swarm algorithm and can also achieve good convergence results. In fact, $v_i^{k+1}$ is not only influenced by these three quantities.
The position ($x_j^k$, $j\neq i$) of other individuals can also have an impact on it. Although the impact is weak, we can still not ignore these information.
How to embed the position information of all other particles ($x_j^k$, $j\neq i$)  into each iteration process (Eqs.(\ref{iteratevelocity}) and (\ref{iterateposition}))? And how to ensure that the position and velocity information of the $i$-th particle is influenced by all other particles?How to design such a new particle swarm algorithm is the problem that this article needs to consider. Based on a series of three experiments, show that in most cases this new particle swarm algorithm will have better optimization performance than the PSO.

The particle swarm optimization algorithm can be understood as a spatiotemporal system. $i$ and $k$ represent space and time, respectively. Therefore, can it be connected to another famous spatiotemporal system - the coupled map lattice system? The answer is affirmative.

The important way to improve particle swarm optimization is to modify its neighborhood topology structure \cite{Suganthan,Mendesfully,Kennedyfully,Lukasikfully,KennedyMenders,KennedyJamesSmall,ReyesMedina,Liu,LiuQunfeng,Nguyen},
These topological structures include but are not limited to the following types:  Small-world network, ring topology, fully connected topology, star topology, mesh Topology,
tree topology, toroidal topology, flustering algorithms, dynamic neighborhood topology, and so on.

For example, the fully informed PSO adopts the following velocity update formula, while the position update formula remains unchanged \cite{Suganthan,Mendesfully,Kennedyfully,Lukasikfully}.
\begin{equation}
\textbf{\textit{v}}_i\leftarrow \chi \left(\textbf{\textit{v}}_i+\sum\limits_{n=1}^{N_i}\frac{U(0,\varphi)(\textbf{\textit{p}}_{nbr(n)}-\textbf{\textit{x}}_i)}{N_i}\right).  \label{fullyPSOVec}
\end{equation}
where $\chi$ known as a constriction factor. $v_i$ is the velocity vector for individual $i$.
$U(0, \varphi)$ is a uniform random number generator.
$N_i$ is the number of neighbors particle $i$ has and $nbr(n)$ is $i's$ $nth$ neighbor. If $nbr(n)$ includes only $i$ itself and its
best neighbor, then this formula degenerates to the PSO.
$\textbf{\textit{p}}_{nbr(n)}$ is the best solution in the $nth$ neighbourhood of the particle $\textbf{\textit{x}}_i$.

In \cite{stereotyping,LiangCluster,Yazawa,GaoYing}, clustering algorithm is integrated into PSO.
For example, the speed update formula is as follows \cite{GaoYing}.
\begin{equation}
  v_i^{(t+1)}=v_i^{(t)}+c_0r_0(p_i^{(t)}-x_i^{(t)})+\sum_{j=1}^Kc_jr_j(q_j^{(t)}-x_i^{(t)}).\label{Gaocluster}
\end{equation}
where $K$ represents the number of subgroups in the clustering process.
$q_j^{(t)}$ means the optimal solution in the $j$-th subgroup during the $t$-th iteration process.
During each iteration, the particle swarm is divided into several subgroups using clustering methods, and the optimal solution in each subgroup is integrated into the PSO.
The performance of particle swarm optimization based on clustering method is superior to that of PSO.

In \cite{allparticle}, the speed update formula is as follows.
\begin{equation}
 v^k(t+1)=\omega\times v^k(t)+c\times \sum_{i=1}^N \lambda_i\times (p_i^k(t)-x^k(t)).\label{allvelocity}
\end{equation}
where $N$ is sarm size. In the first $t$ iterations, $p_i^k(t)$ represents the current optimal solution value of the $k$-th dimensional component of the $i$-th particle.
During each iteration, the local optimal solutions of all particles participate in the velocity update.

\section{Coupled Map Lattice}

A coupled map lattice(CML) is a high-dimensional discrete dynamical system formed by a given discrete dynamical system (mapping) interacting (coupling) in a certain way.
It was first proposed by Kuniko Kaneko\cite{Kaneko}.

The first simple CML is proposed for studying spatiotemporal chaos: considering the phenomena generated by local chaotic processes and spatial diffusion processes. Take one-dimensional mapping as the simplest representative chaos and use discrete Laplace operator for diffusion. Some phenomena in nature can be represented by reaction-diffusion equations, namely
\begin{equation}\label{Reactdiff}
\partial_t u=F(u)+\epsilon\nabla^2 u
\end{equation}
The local reaction process can be represented by a nonlinear mapping $x(i)\rightarrow x'(i)=f(x(i))$, for example, Logistic mapping $f(x)=rx(1-x)$.
And the diffusion process is obtained by the discrete Laplace operator
$$x'(i)\rightarrow (1-\epsilon)x'(i)+\frac{\epsilon}{2}\left[x'(i+1)+x'(i-1)\right].$$
Combine the above local reaction processes with the diffusion process represented by the discrete Laplace operator to obtain a coupled map lattice
\begin{equation}\label{basicCML}
x_{n+1}(i)=(1-\epsilon)f(x_n(i))+\frac{\epsilon}{2}\left[f(x_n(i-1))+f(x_n(i+1))\right]
\end{equation}
$n$ represents discrete time. $x$ is system state. $i$ is lattice coordinates ($i=1, 2, \cdots, L$). $L$ is system size. $\epsilon$ is coupling parameter (coupling strength), and $f$ is a local function or lattice function\cite{Kaneko,YuanCML}.
Usually, periodic boundary conditions are taken, i.e., $x_n(0)=x_n(L).$

Coupled map lattice with different types of topological structures have been proposed.
The selection of coupling topology structure will affect the dynamics of the entire system.
Coupled map lattice is relatively simple, but it has most of the characteristics of spatiotemporal chaos \cite{Kaneko,YuanCML,Vasegh,Khellat,Cosenza,travelingwaves,Neufeld,Spatiotemporal,Leonid}.
Usually, the calculations for each lattice point are exactly the same, and some results of the theory of low dimensional dynamical systems can still be applied to their analysis.
Coupled map lattices and their applications have been widely studied.
Fore example, some common topological structures of couple map lattices are as follows.

Globally coupled map lattice \cite{Khellat,Cosenza,Popovych}
\begin{equation}\label{GlobalCML}
x_{n+1}(i)=(1-\epsilon)f(x_n(i))+\frac{\epsilon}{N-1}\sum_{j=1,j\neq i}^N f(x_n(j)).
\end{equation}
The information at time $n+1$ of the $i$-th lattice point is partially influenced by the information of all lattice points at time $n$, $1-\epsilon$ and $\frac{\epsilon}{N-1}$ can be understood as weight coefficient.

Accumulated coupled map lattice \cite{Vasegh}
\begin{align}
x_{n+1}(1)&=f(x_n(1)),\label{AccumulatedCML1}\\
x_{n+1}(j)&=(1-\epsilon)f(x_n(j))+\frac{\epsilon}{j-1}\sum_{i=1}^{j-1} f(x_n(i))\label{AccumulatedCML2}.
\end{align}

The information at time $n+1$ of the $j$-th lattice point is partially influenced by the information of the first $j$ lattice points and itself at time $n$. It is the accumulation of information from the previous $j$ lattice points.

There are many forms of coupled map lattice systems, such as cross coupled map lattice systems, two-dimensional coupled map lattice systems, etc.
In coupled map lattice systems, the adopted lattice functions are usually one-dimensional chaotic maps, such as logistic maps, Tent maps, etc.
Therefore, coupled map lattice systems have complex dynamic properties and it has a wide rang of applications, especially in the field of chaotic image encryption.
The mutual influence of information in coupled map lattice systems inspires us to apply this idea to particle swarm optimization algorithms and integrate it with the foraging process of birds. We link the information exchange process of birds to a coupled map lattice systems.

\section{Globally Coupled Particle Swarm Optimization}
Particle swarm optimization algorithm is a method of information sharing among the birds, and continuously updating information of the birds by searching for the optimal bird's position information.
Coupled map lattice system is the mutual influence between lattice points, which is somewhat similar to the mutual influence between birds.
Therefore, we attempt to integrate the coupled map lattice with the particle swarm optimization algorithm to propose a new modified particle swarm optimization algorithm.
Now we will embed the concept of globally coupled map lattice $(\ref{GlobalCML})$ into particle swarm optimization algorithm $(\ref{iteratevelocity})$ and $(\ref{iterateposition})$.
A new particle swarm optimization algorithm will be proposed. We name it particle swarm optimization with coupled map lattice (denoted as PSOCML).

\begin{figure}[h]
\centering
\includegraphics[width=0.6\textwidth]{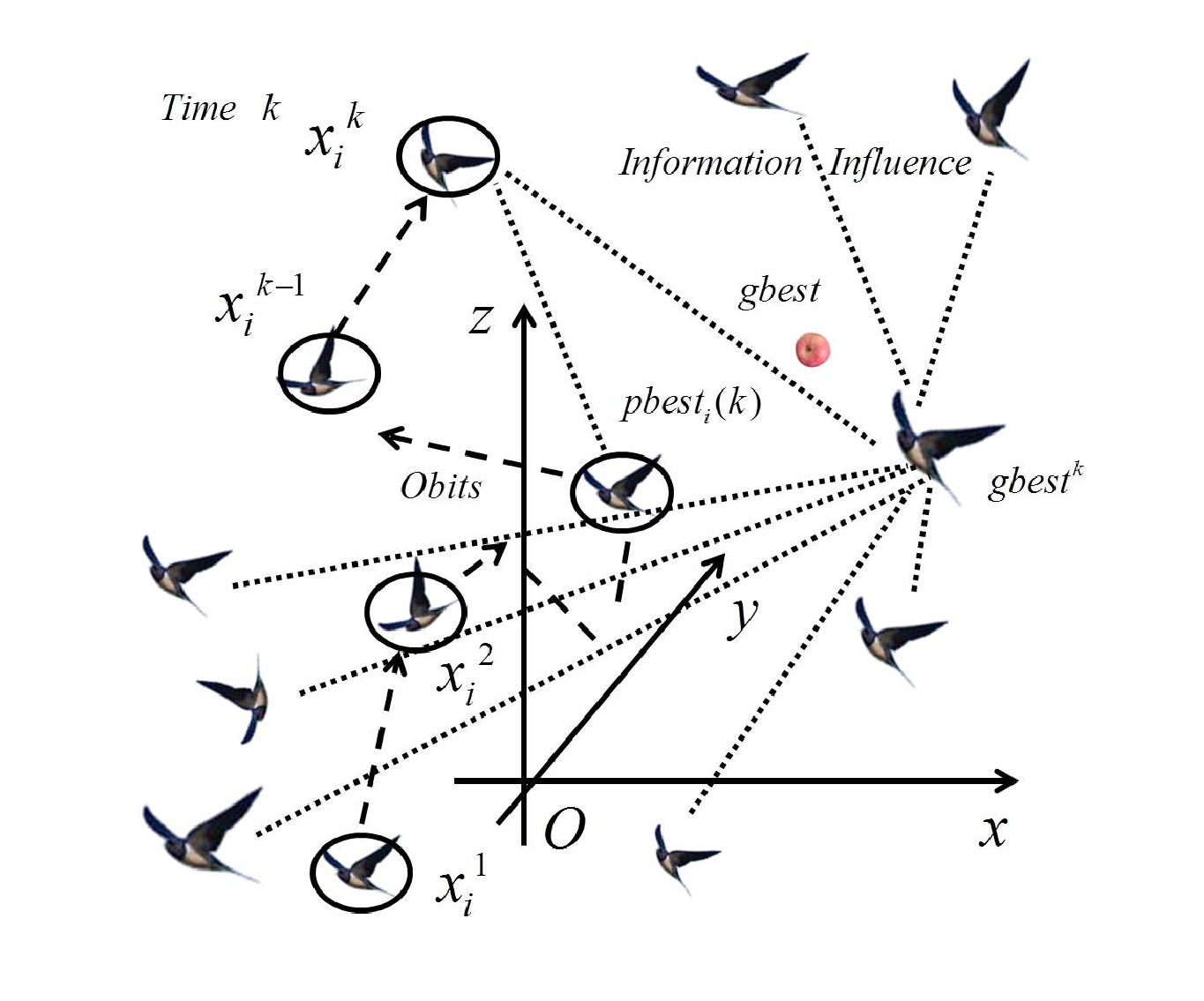}
\caption{Equation (\ref{iterateposition}) of position update }
\label{velocityadd}
\end{figure}

Two update iteration formulas of PSOCML are as shown in equations $(\ref{PSOwithGCML})$ and $(\ref{iteratepositionwithGCML})$,
where the velocity update formula $(\ref{PSOwithGCML})$ undergoes significant changes.
The coupling principle of the globally coupled map lattice $(\ref{GlobalCML})$ is utilized here.
\begin{align}
v_i^{k+1}&=w\times v_i^k+c_1 r_1(pbest_i^k-x_i^k)\nonumber\\
&~~~+c_2 r_2 \left((1-\epsilon)(gbest^k-x_i^k)+\frac{\epsilon}{n-1}\sum\limits_{j=1,j\neq i}^n (gbest^k-x_j^k)\right), \label{PSOwithGCML}\\
x_i^{k+1}&=x_i^k+v_i^{k+1},~~~i=1,2,\cdots,n,\label{iteratepositionwithGCML}
\end{align}
where $n$ is the swarm size, $v_i^{k+1}$ and $x_i^{k+1}$ represent the velocity and position of particle $i$ at $(k-1)th$ iteration step respectively.
$\epsilon$ is the coupling strength. $(1-\epsilon)(gbest^k-x_i^k)$ indicates that information of the global optimal position $gbest^k$  affects the current particle $x_i^k$,
and $\frac{\epsilon}{n-1}\sum\limits_{j=1,j\neq i}^n (gbest^k-x_j^k)$ describes that the global optimal position has an impact on other particles.
$\left((1-\epsilon)(gbest^k-x_i^k)+\frac{\epsilon}{n-1}\sum\limits_{j=1,j\neq i}^n (gbest^k-x_j^k)\right)$ is social part.
Compared to particle swarm optimization, it contains  much information of swarm. Denote $c_1r_1=k_1$, $c_2r_2(1-\epsilon)=k_2$ and $\frac{c_2r_2\epsilon}{n-1}=k_3$,
the schematic diagram of equation (\ref{PSOwithGCML}) is shown in Fig.\ref{velocityadd}. These also demonstrate that everything is interconnected, and each bird has an influence on other birds, and it is also influenced by others. Although the form of (\ref{PSOwithGCML}) is more complex than equation $(\ref{iteratevelocity})$ in particle swarm optimization algorithm, the principle is still simple and effective. This algorithm will have better ability to find the optimal value.
$\frac{\epsilon}{n-1}\sum\limits_{j=1,j\neq i}^n (gbest^k-x_j^k)$ represents all other particles $x_j^k$ ($j\neq i$) are drawn towards $gbest^k$.
These information also affects the flight speed $v_i^{k+1}$ of particle $x_i$. $\epsilon$ ($0\leq \epsilon\leq 1$) is the weight value of information, usually $\epsilon$ has a smaller value.

The pseudocode PSOCML is explained in Algorithm.
\begin{table}[ht]
  \centering
  \begin{tabular}{p{13cm}}
    \toprule
    {\bf Algorithm} Globally Coupled Particle Swarm Optimization\\
    \midrule
  {\bf Initialization} of the population position $x_i^k$, velocities $v_i^k$ and parameters $c_1$, $c_2$,
  $n$, $w$, $r_1$, $r_2$, $\epsilon$, etc; \\
  Compute the objective function value (\ref{miniproblem}) of each initialization position $x_i^k$;\\
  Identify the global best value ($gbest^k$), and previous best value ($pbest_i^k$) of each \\
  particle, where $pbest_i^k=x_i^k$, $k=1$ and $i=1, 2, \cdots, n$;\\
   {\bf while} termination criteria (Maximum number of iterations) {\bf do} \\
   ~~~~{\bf for} each particle $x_i^k$ and velocity $v_i^k$ {\bf do} \\
   ~~~~~~(i) Based on equation (\ref{PSOwithGCML}), update speed velocity $v_i^{k+1}$;\\
   ~~~~~~(ii) Based on equation (\ref{iteratepositionwithGCML}), update speed velocity $x_i^{k+1}$;\\
   ~~~~{\bf end for}\\
   ~~~~With updated solutions $x_i^{k+1}$, calculate the objective function value;\\
   ~~~~Update the positions of global and individual optimal solutions once again,\\
   ~~~~i.e., $gbest^k$ and $pbest_i^k$\\
   {\bf end while}\\
  Return the optimal value $gbest^k$ as the optimal solution of the optimization problem (\ref{miniproblem}).\\
  {\bf Until} stopping condition is true.\\
    \bottomrule
  \end{tabular}
\end{table}

\begin{remark}\label{specialone}
When $\epsilon=0$,  the equations $(\ref{PSOwithGCML})$ and $(\ref{iteratepositionwithGCML})$ degenerate into PSO.
\end{remark}

\begin{remark}\label{reducetopso}
When $\epsilon=1$, the equation $(\ref{PSOwithGCML})$ becomes
\begin{equation}\label{specialcase2}
v_i^{k+1}=w\times v_i^k+c_1 r_1(pbest_i^k-x_i^k)
+\frac{c_2 r_2}{n-1}\sum\limits_{j=1,j\neq i}^n (gbest^k-x_j^k).
\end{equation}
This situation (\ref{specialcase2}) will be analyzed in the future. Whether it is superior to PSO is still unknown.
\end{remark}

\begin{remark}\label{psoCMLandworkerantlaw}
Coupled map lattice exhibits spatiotemporal chaos \cite{Kaneko,YuanCML,Vasegh,Khellat,Cosenza,travelingwaves,Neufeld,Spatiotemporal,Leonid}, and its orbit has ergodicity.
When the time is long enough, the orbit will fill the entire phase space.
These chaotic features are suitable for updating particle communities in PSO \cite{psochaos2005}.
kanta Matsumoto and Chihiro Ikuta proposed the new PSO with CML and worker ant's law \cite{Matsumoto}.
The CML is applied to the moving equation of PSO\cite{Matsumoto,KantaMChihiroI}, and the modification of the velocity update formula is as follows.
\begin{align}
v_{id}(t+1)=&\omega(t)v_{id}(t)+C_1 rand()[x_{pbest}(t)-x_{id}(t)]\nonumber\\
&+C_2 rand()[x_{gbest}(t)-x_{id}(t)]+\eta(t)f[z_{id}(t)].
\end{align}
Where $\eta(t)$ is the amplitude controller of output of CML.
$f[z_{id}(t)]$ is chaotic orbits of CML.$\eta(t)f[z_{id}(t)]$
is independent from $x_{pbest}(t)$, $x_{gbest}(t)$ and each particle $x_{id}(t)$.
Adding this term ($\eta(t)f[z_{id}(t)]$) has disturbance effect on the velocity value
to avoid the algorithm getting stuck in local optima too early.
Here CML is added to the velocity equation of each particle. CML helps PSO escaping out from the local solutions.
CML enhances the search capability of the PSO.
This method involves chaotic perturbation of velocity. This increases ergodicity of speed,
thereby affecting the positional information of particles.
This is different from various topological structures of PSO \cite{Suganthan,Mendesfully,Kennedyfully,Lukasikfully,KennedyMenders,KennedyJamesSmall,ReyesMedina,Liu,LiuQunfeng,Nguyen}.
\end{remark}

\section{Experiments and Results}
The comparison is made with the PSO, and other variations of PSO.
PSO and most of its improved algorithms are designed to solve single objective global optimal problems.
These can not find multiple optimal points of a multimodal function.
Therefore, the ones we have selected below are all single objective testing functions.

The performance of the GCPSO is tested on modern benchmark functions, i.e., CEC 2017 (Awad et al., 2016). The performance
of the GCPSO is compared with PSO and other improved PSO.

In the second test suite ($F1-F30$), CEC2017 test suite (Awad et al.,
2016) (F2 has been excluded because it shows unstable behavior especially for higher dimensions, and significant performance variations
for the same algorithm implemented in Matlab) are selected to evaluate
the performance the proposed GCPSO. 
Originally CEC $2017$ consisted of $30$ functions. Then, function number 2 was excluded from the competition. At the time of writing this note (March 2021),
an updated problem definitions renumbered functions formerly used in $2017$ (old function $3$ is now called $2$, etc.).
$https://staff.elka.pw.edu.pl/~rbiedrzy/publ/IPOP.html$
These benchmark functions are
classified into four categories: unimodal functions ($F1-F3$), multimodal
functions ($F4-F10$), hybrid functions ($F11-F20$), composition functions
($F21-F30$). These functions are used to further verify the validity of the
proposed GCPSO. The problem dimension, population size
and max number of fitness evaluations are set uniformly as 30, 40 and
$3.0\times 105$. All the PSO algorithms are run for $51$ independent times.
These shifted and rotated functions are more complex
and make our test results more convincing.

\begin{figure}
	\centering
	\includegraphics[width=7cm]{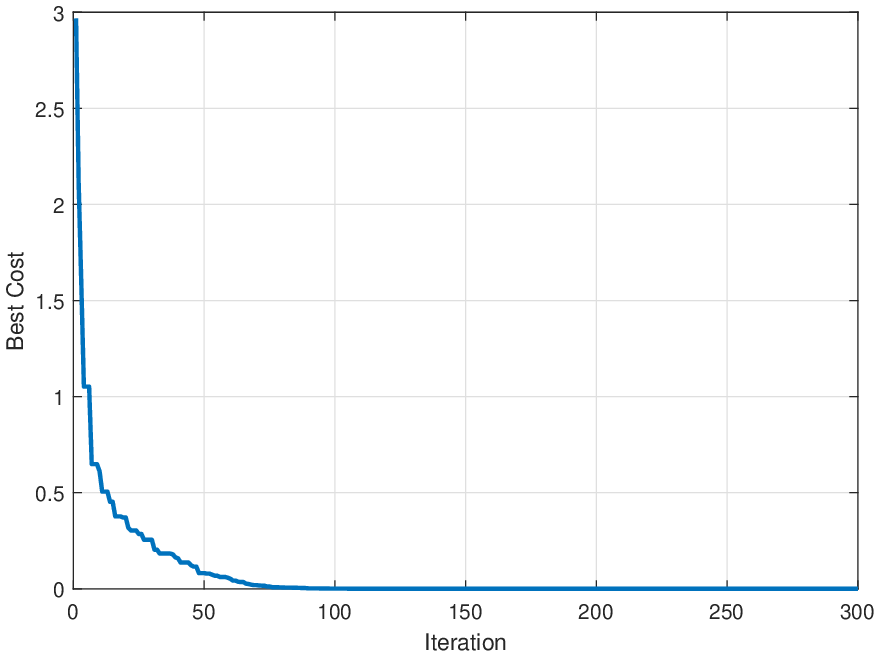}	
	\caption{Ackley test function,Iteration $300$: Best Cost$=1.4655e-14$}
	\label{Ackley}
\end{figure}

\begin{figure}
	\centering
	\includegraphics[width=7cm]{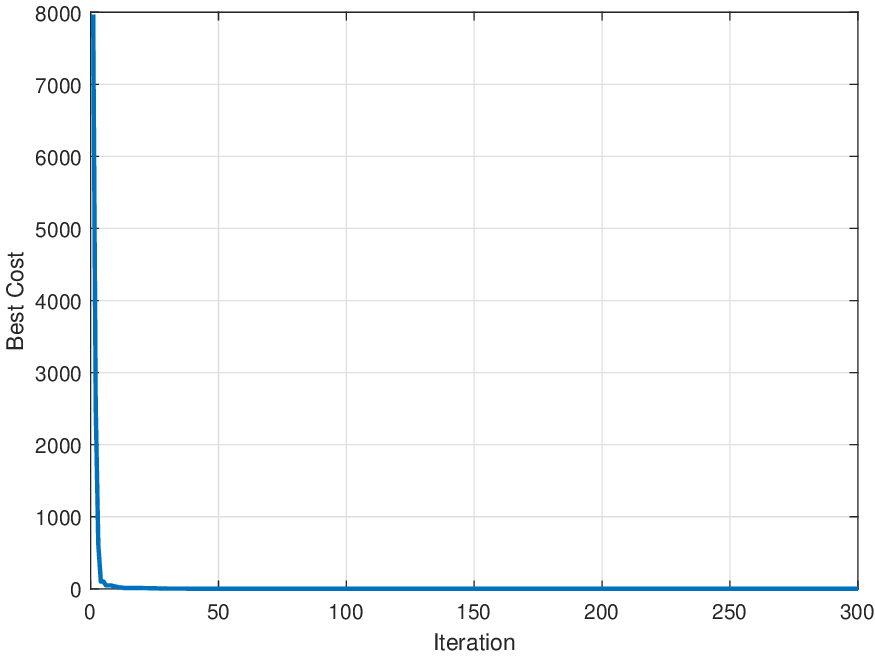}	
	\caption{DeJongf4 test function, Iteration $300$: Best Cost $= 5.1349e-49$}
	\label{Ackley}
\end{figure}

\section{Conclusion}\label{conclsec}
We integrate the framework of the globally coupled map lattice system into the particle swarm optimization algorithm and propose a new set of globally coupled particle swarm optimization algorithms. After experimental testing, the new globally coupled particle swarm optimization  has excellent optimization ability.

\subsection*{Acknowledgements}
This work is supported by  (No. ).


\begin{thebibliography}{99}
\setlength{\itemsep}{-2pt}
\scriptsize



\bibitem{Kennedy1995}Kennedy J, Eberhart R. Particle swarm optimization[C]//Proceedings of ICNN'95-international conference on neural networks. ieee, 1995, 4: 1942-1948.
\bibitem{Eberhast1995}Eberhart R, Kennedy J. A new optimizer using particle swarm theory[C]//MHS'95. Proceedings of the sixth international symposium on micro machine and human science. Ieee, 1995: 39-43.




\bibitem{parallelPSO}Charilogis V, Tsoulos IG, Tzallas A. An improved parallel particle swarm optimization. SN Computer Science. 2023, 4(6):766.
\bibitem{30objectivefunctions}Plevris V, Solorzano G. A Collection of 30 Multidimensional Functions for Global Optimization Benchmarking. Data. 2022; 7(4):46.
\bibitem{Anula}Anula Khare, Saroj Rangnekar. A review of particle swarm optimization and its applications in solar photovoltaic system. Applied Soft Computing 13.5 (2013): 2997-3006.
\bibitem{Houssein}Houssein, E. H., Gad, A. G., Hussain, K., $\&$ Suganthan, P. N. (2021). Major advances in particle swarm optimization: theory, analysis, and application. Swarm and Evolutionary Computation, 63, 100868.
\bibitem{AlRashidi}M. R. AlRashidi, M. E. El-Hawary, A Survey of Particle Swarm Optimization Applications in Electric Power Systems, in IEEE Transactions on Evolutionary Computation, vol. 13, no. 4, pp. 913-918, Aug. 2009.
\bibitem{Shi}Shi, Yuhui. Particle swarm optimization: developments, applications and resources. Proceedings of the 2001 congress on evolutionary computation (IEEE Cat. No. 01TH8546). Vol. 1. IEEE, 2001.
\bibitem{Kulkarni}Kulkarni, Mr Ninad K., et al. Particle swarm optimization applications to mechanical engineering-A review. Materials Today: Proceedings 2.4-5 (2015): 2631-2639.
\bibitem{Ekrem}Ekrem \"{o}zge, Bekir Aksoy. Trajectory planning for a 6-axis robotic arm with particle swarm optimization algorithm. Engineering Applications of Artificial Intelligence 122 (2023): 106099.
\bibitem{Ahmed}Ahmed G. Gad. Particle swarm optimization algorithm and its applications: a systematic review. Archives of computational methods in engineering 29.5 (2022): 2531-2561.
\bibitem{Tareq}Tareq M. Shami, Ayman A. El-Saleh, et al. Particle swarm optimization: A comprehensive survey. Ieee Access 10 (2022): 10031-10061.
\bibitem{Tiwari}Tiwari, Sukriti, and Ashwani Kumar. Advances and bibliographic analysis of particle swarm optimization applications in electrical power system: concepts and variants. Evolutionary Intelligence 16.1 (2023): 23-47.
\bibitem{Nayak}Nayak, Janmenjoy, et al. 25 years of particle swarm optimization: Flourishing voyage of two decades. Archives of Computational Methods in Engineering 30.3 (2023): 1663-1725.
\bibitem{Ramirez}Ramirez-Figueroa, John A., et al. A new principal component analysis by particle swarm optimization with an environmental application for data science. Stochastic Environmental Research and Risk Assessment 35.10 (2021): 1969-1984.
\bibitem{Xue}Xue Yu, Qi Zhang, Adam Slowik. Automatic topology optimization of echo state network based on particle swarm optimization. Engineering Applications of Artificial Intelligence 117 (2023): 105574.
\bibitem{Pace}Pace Francesca, Alessandro Santilano, Alberto Godio. A review of geophysical modeling based on particle swarm optimization. Surveys in Geophysics 42.3 (2021): 505-549.
\bibitem{Kashani}Kashani, Ali R., et al. Particle swarm optimization variants for solving geotechnical problems: review and comparative analysis. Archives of Computational Methods in Engineering 28 (2021): 1871-1927.
\bibitem{Sedighizadeh}Sedighizadeh Davoud, et al. GEPSO: A new generalized particle swarm optimization algorithm. Mathematics and Computers in Simulation 179 (2021): 194-212.
\bibitem{Felix}Felix T.S. Chan, Manoj Tiwari, eds. Swarm Intelligence: focus on ant and particle swarm optimization. BoD$-$Books on Demand, 2007.
\bibitem{Mercang}Burcu Adiguzel Mercangoz. Applying particle swarm optimization: new solutions and cases for optimized portfolios. Vol. 306. Springer Nature, 2021.
\bibitem{Micael}Micael Couceiro , Pedram Ghamisi. Fractional Order Darwinian Particle Swarm Optimization Applications and Evaluation of an Evolutionary Algorithm, Springer, 2016.
\bibitem{Olsson}Olsson, Andrea E. Particle swarm optimization: theory, techniques and applications. Nova Science Publishers, Inc., 2010.
\bibitem{Mikki}Mikki  Said M., Ahmed A. Kishk. Particle swarm optimizaton: a physics-based approach. Springer Nature, 2022.
\bibitem{Eberhart}Eberhart, Russell C., Yuhui Shi, and James Kennedy. Swarm intelligence. Elsevier, 2001.
\bibitem{Sengupta}Sengupta S, Basak S, Peters R A. Particle Swarm Optimization: A survey of historical and recent developments with hybridization perspectives. Machine Learning and Knowledge Extraction, 2018, 1(1): 157-191.


\bibitem{Mendesfully}Mendes, Rui, James Kennedy, and Jos$\acute{e}$ Neves. The fully informed particle swarm: simpler, maybe better. IEEE transactions on evolutionary computation 8.3 (2004): 204-210.
\bibitem{Kennedyfully}Kennedy, James, and Rui Mendes. Neighborhood topologies in fully informed and best-of-neighborhood particle swarms. IEEE Transactions on Systems, Man, and Cybernetics, Part C (Applications and Reviews) 36.4 (2006): 515-519.
\bibitem{Lukasikfully}Lukasik, Szymon, and Piotr A. Kowalski. Fully informed swarm optimization algorithms: basic concepts, variants and experimental evaluation. 2014 Federated Conference on Computer Science and Information Systems. IEEE, 2014.
\bibitem{KennedyMenders}Kennedy, James, and Rui Mendes. Population structure and particle swarm performance. Proceedings of the 2002 Congress on Evolutionary Computation. CEC'02 (Cat. No. 02TH8600). Vol. 2. IEEE, 2002.
\bibitem{KennedyJamesSmall}Kennedy, James. Small worlds and mega-minds: effects of neighborhood topology on particle swarm performance. Proceedings of the 1999 congress on evolutionary computation-CEC99 (Cat. No. 99TH8406). Vol. 3. IEEE, 1999.
\bibitem{ReyesMedina}Reyes-Medina, Angelina Jane, Gregorio Toscano Pulido, and Jos$\acute{e}$ Gabriel Ram$\acute{i}$rez-Torres. A comparative study of neighborhood topologies for particle swarm optimizers. International conference on evolutionary computation. Vol. 2. SciTePress, 2009.
\bibitem{Liu}Liu Y, Zhao Q, Shao Z, et al. Particle swarm optimizer based on dynamic neighborhood topology[C]//Emerging Intelligent Computing Technology and Applications. With Aspects of Artificial Intelligence: 5th International Conference on Intelligent Computing, ICIC 2009 Ulsan, South Korea, September 16-19, 2009 Proceedings 5. Springer Berlin Heidelberg, 2009: 794-803.
\bibitem{LiuQunfeng}Liu, Qunfeng, et al. Topology selection for particle swarm optimization. Information Sciences 363 (2016): 154-173.
\bibitem{Nguyen}Nguyen, Loc, Ali A. Amer, Hassan I. Abdalla. A general framework of particle swarm optimization. Proceedings of the Future Technologies Conference. Cham: Springer International Publishing, 2022.
\bibitem{Suganthan}Suganthan, Ponnuthurai N. Particle swarm optimiser with neighbourhood operator. Proceedings of the 1999 congress on evolutionary computation-CEC99 (Cat. No. 99TH8406). Vol. 3. IEEE, 1999.
\bibitem{allparticle}Liu, Qing, Jin Li, Haipeng Ren, and Wei Pang. All particles driving particle swarm optimization: Superior particles pulling plus inferior particles pushing. Knowledge-Based Systems 249 (2022): 108849.

\bibitem{stereotyping}Kennedy J. (2000,July). Stereotyping: Improving particle swarm performance with cluster analysis. In Proceedings of the 2000 congress on evolutionary computation. CEC00 (Cat. No. 00TH8512) (Vol. 2, pp. 1507-1512). IEEE.
\bibitem{LiangCluster}Liang, Xiaolei, et al. An adaptive particle swarm optimization method based on clustering. Soft Computing 19 (2015): 431-448.
\bibitem{Yazawa}Yazawa, Kazuyuki, Makoto Motoki, Keiichiro Yasuda. Cluster-structured particle swarm optimization with interaction. 2009 ICCAS-SICE. IEEE, 2009.
\bibitem{GaoYing}Ying Gao, Shengli Xie, Ruoning Xu, et al. A multi sub-population particle swarm optimizer based on clustering (in chinese), Application Research of Computers, 4(2006):40-41.

\bibitem{psochaos2005}Bo Liu, Ling Wang,Yi-Hui Jin, et al. Improved particle swarm optimization combined with chaos. Chaos, Solitons $\&$ Fractals, 2005, 25(5): 1261-1271.
\bibitem{Matsumoto}Matsumoto, Kanta, Chihiro Ikuta. PSO with Coupled Map Lattice and Worker Ant's Law. 2021 IEEE Congress on Evolutionary Computation (CEC). IEEE, 2021.
\bibitem{KantaMChihiroI}Kanta M, Chihiro I. Particle Swarm Optimization with Coupled Map Lattice. IEICE Proceeding Series, 2020, 74: 389-392.



\bibitem{Kaneko}K. Kaneko, Theory and Applications of Coupled Map Lattices, Nonlinear Science: Theory and Applications, JohnWiley $\&$ Sons, New York, NY, USA, 1993.
\bibitem{YuanCML}Liguo Yuan, Qigui Yang. A proof for the existence of chaos in diffusively coupled map lattices with open boundary conditions. Discrete Dynamics in Nature and Society, 2011,
Article ID 174376, 16 pages.
\bibitem{Vasegh}Vasegh Nastaran. Spatiotemporal and synchronous chaos in accumulated coupled map lattice. Nonlinear Dynamics 89.2 (2017): 1089-1097.
\bibitem{Khellat}Khellat Farhad, Akashe Ghaderi,  Nastaran Vasegh. Li-Yorke chaos and synchronous chaos in a globally nonlocal coupled map lattice. Chaos, Solitons $\&$ Fractals 44.11 (2011): 934-939.
\bibitem{Cosenza}Cosenza, M. G., and A. Parravano. Turbulence in globally coupled maps. Physical Review E 53.6 (1996): 6032.
\bibitem{Popovych}Popovych O., Yu Maistrenko, Erik Mosekilde. Loss of coherence in a system of globally coupled maps. Physical Review E 64.2 (2001): 026205.
\bibitem{travelingwaves}Kaneko Kunihiko. Chaotic traveling waves in a coupled map lattice. Physica D: Nonlinear Phenomena 68.3-4 (1993): 299-317.
\bibitem{Neufeld}Neufeld Z., T. Vicsek. Spatiotemporal chaos in a coupled map lattice with unstable couplings. Journal of Physics A: Mathematical and General 28.18 (1995): 5257.
\bibitem{Spatiotemporal}Kaneko Kunihiko. Spatiotemporal chaos in one-and two-dimensional coupled map lattices. Physica D: Nonlinear Phenomena 37.1-3 (1989): 60-82.
\bibitem{Leonid}Bunimovich Leonid A.,  Ya G. Sinai. Spacetime chaos in coupled map lattices. Nonlinearity 1.4 (1988): 491.


\bibitem{testShakti}Chourasia S, Sharma H, Singh M, et al. Global and local neighborhood based particle swarm optimization[C]//Harmony Search and Nature Inspired Optimization Algorithms: Theory and Applications, ICHSA 2018. Springer Singapore, 2019: 449-460.
\bibitem{test30functions} Vagelis Plevris, German Solorzano. A collection of $30$ multidimensional functions for global optimization benchmarking. Data, 2022: 46.
\bibitem{CEC2005}Suganthan, P. N., Hansen, N., Liang, J. J., Deb, K., Chen, Y.-P., Auger,A., Tiwari, S., 2005. Problem definitions and evaluation criteria for the
CEC 2005 special session on real-parameter optimization. Tech. Rep. 2005005, Nanyang Technological University and KanGAL, Singapore and IIT Kanpur, India.








\bibitem{PSOKennedy}Kennedy, James, and Russell Eberhart. Particle swarm optimization. In Proceedings of ICNN'95-international conference on neural networks, vol. 4, pp. 1942-1948. IEEE, 1995.
\bibitem{ShiRussell}Shi, Yuhui, and Russell Eberhart. A modified particle swarm optimizer. In 1998 IEEE international conference on evolutionary computation proceedings. IEEE world congress on computational intelligence (Cat. No. 98TH8360), pp. 69-73. IEEE, 1998.
\bibitem{Bratton}Bratton, Daniel, and James Kennedy. Defining a standard for particle swarm optimization. 2007 IEEE swarm intelligence symposium. IEEE, 2007.
\bibitem{Poli}Poli, Riccardo, James Kennedy, and Tim Blackwell. Particle swarm optimization: An overview. Swarm intelligence 1(2007): 33-57.
\bibitem{Vanneschi}Vanneschi, Leonardo, and Sara Silva. Particle Swarm Optimization. Lectures on Intelligent Systems. Cham: Springer International Publishing, 2023. 105-111.
\bibitem{Nayak}Nayak, Janmenjoy, et al. 25 years of particle swarm optimization: Flourishing voyage of two decades. Archives of Computational Methods in Engineering 30.3 (2023): 1663-1725.
\bibitem{Chih}Chih, Mingchang. Stochastic stability analysis of particle swarm optimization with pseudo random number assignment strategy. European Journal of Operational Research 305.2 (2023): 562-593.
\bibitem{ClercPSO}Clerc, Maurice. Particle swarm optimization. Vol. 93. John Wiley \& Sons, 2010.
\bibitem{Shi}Shi, Yuhui. Particle swarm optimization. IEEE connections 2.1 (2004): 8-13.
\bibitem{overviewpso}Wang, Dongshu, Dapei Tan, and Lei Liu. Particle swarm optimization algorithm: an overview. Soft computing 22 (2018): 387-408.
\bibitem{comprehensive}Shami, Tareq M., et al. Particle swarm optimization: A comprehensive survey. IEEE Access 10 (2022): 10031-10061.
\bibitem{Gad}Gad, Ahmed G. Particle swarm optimization algorithm and its applications: a systematic review. Archives of computational methods in engineering 29.5 (2022): 2531-2561.
\bibitem{Yuan}Liguo Yuan, Qigui Yang. Parameter identification and synchronization of fractional-order chaotic systems. Communications in Nonlinear Science and Numerical Simulation 17.1 (2012): 305-316.



\end{thebibliography}
\end{document}